\documentclass{amsart}

\newtheorem{theorem}{Theorem}[section]
\newtheorem{lemma}[theorem]{Lemma}
\newtheorem{proposition}[theorem]{Proposition}
\newtheorem{corollary}[theorem]{Corollary}
\newtheorem*{claim}{Claim}

\theoremstyle{definition}
\newtheorem{definition}[theorem]{Definition}
\newtheorem{example}[theorem]{Example}

\newtheorem{remark}[theorem]{Remark}

\newtheorem{question}[theorem]{Question}
\newtheorem*{acknowledgement}{Acknowledgement}

\theoremstyle{remark}

\newcommand{\NN}{\mathbb{N}}
\newcommand{\ZZ}{\mathbb{Z}}
\newcommand{\QQ}{\mathbb{Q}}
\newcommand{\RR}{\mathbb{R}}

\newcommand{\PP}{\mathbb{P}}
\renewcommand{\AA}{\mathbb{A}}

\newcommand  {\shA}     {\mathcal{A}}

\newcommand  {\shB}     {\mathcal{B}}

\newcommand  {\shE}     {\mathcal{E}}
\newcommand  {\shF}     {\mathcal{F}}
\newcommand  {\shG}     {\mathcal{G}}

\newcommand  {\shI}     {\mathcal{I}}

\newcommand  {\shK}     {\mathcal{K}}
\newcommand  {\shM}     {\mathcal{M}}

\newcommand  {\shL}     {\mathcal{L}}

\newcommand  {\fob}     {\mathfrak{b}}

\newcommand  {\fom}     {\mathfrak{m}}

\newcommand  {\fop}     {\mathfrak{p}}

\newcommand  {\foX}     {\mathfrak{X}}

\newcommand  {\aff}     {{\text{aff}}}

\newcommand  {\Card}    {\operatorname{Card}}

\newcommand  {\Hom}     {\operatorname{Hom}}

\newcommand  {\limdir}  {\varinjlim}

\newcommand  {\lra}     {\longrightarrow}

\renewcommand{\O}       {\mathcal{O}}

\newcommand  {\ord}     {\operatorname{ord}}

\newcommand  {\pr}      {\operatorname{pr}}
\newcommand  {\Proj}    {\operatorname{Proj}}

\newcommand  {\quadand} {\quad\text{and}\quad}
\newcommand  {\Quot}    {\operatorname{Quot}}

\newcommand  {\ra}      {\rightarrow}

\newcommand  {\rank}    {\operatorname{rank}}

\newcommand  {\Spec}    {\operatorname{Spec}}
\newcommand  {\Spf}     {\operatorname{Spf}}
\newcommand  {\Supp}    {\operatorname{Supp}}

\def\mydate{\number\day\space\ifcase\month \or January\or February\or March\or April\or May\or
June\or July\or August\or September\or October\or November\or
December\fi \space\number\year}

\usepackage{amscd}
\usepackage{amssymb}

\begin{document}

\title[Ample families and homogeneous spectra]
{Ample families, multihomogeneous spectra, and
algebraization of formal schemes}

% Remove or comment out any unused author tags.
% author one information

\author[Holger Brenner]{Holger Brenner}
\address{Department of Pure Mathematics, University of Sheffield, Sheffield S3 7RH, United Kingdom}
%\curraddr{}
\email{h.brenner@shef.ac.uk}

\author[Stefan Schroeer]{Stefan Schr\"oer}
\address{Mathematisches Institut,
Heinrich-Heine-Universit\"at D\"usseldorf,
Universit\"atsstr. 1,
40225 D\"usseldorf, Germany}

\email{schroeer@uni-duesseldorf.de}


\begin{abstract}
Generalizing   homogeneous spectra for
rings graded by  natural numbers, we introduce multihomogeneous spectra for
rings graded by
abelian groups.
Such homogeneous spectra have the same completeness properties
as their classical counterparts, but
are possibly nonseparated. We relate them to  ample families of invertible
sheaves and simplicial toric varieties. As an application, we
generalize Grothendieck's Algebraization Theorem and  show that
formal schemes with certain  ample
families are algebraizable.
\end{abstract}

\maketitle

%===========================================================
\section*{Introduction}

A powerful method to study algebraic varieties is to embed them, if possible,
into some projective space
$\PP^n$. Such an embedding
$X\subset \PP^n$ allows you to view  points
$x\in X$ as homogeneous prime ideals in some
$\NN$-graded ring. The purpose of this paper is to extend this  to
\emph{divisorial varieties}, which are not
necessarily quasiprojective.

The notion of divisorial varieties is due to Borelli \cite{Borelli 1963}.
The class of divisorial varieties
contains all quasiprojective schemes, smooth varieties, and locally
$\QQ$-factorial varieties.
Roughly speaking, divisoriality means that there is a finite collection
of invertible sheaves $\shL_1,\ldots,\shL_r $ so that the whole collection
behaves like
an ample invertible sheaf. Such  collections are  called  \emph{ample
families}.

Our main idea is to define \emph{homogeneous spectra}
$\Proj(S)$ for multigraded rings, that is, for rings graded by an
abelian group of finite type. We obtain
$\Proj(S)$ by patching affine pieces
$D_+(f)=\Spec(S_{(f)})$, where
$f\in S$ are certain homogeneous elements. Roughly speaking, we
demand that  $f$ has many homogeneous divisors.
Multihomogeneous spectra share many properties of classical homogeneous
spectra. For example, they are universally closed, and the
intersection of two affine open subsets is affine. They are,
however, not necessarily separated.

Roberts \cite{Roberts 1998} gave a similar construction
for $\NN^r$-graded rings $S$ satisfying certain
conditions on homogeneous generators.
He used it to study
Hilbert functions in several variables and local multiplicities.
Roberts' homogeneous spectrum is an open subset of ours.

It turns out that a scheme is divisorial if and only if it
admits an embedding into suitable  multihomogeneous
spectra.  More precisely, we shall characterize ample families of
invertible sheaves in terms of
$\Proj(S)$ for the multigraded
ring of global sections
$S=\bigoplus_{d\in\NN^r}
\Gamma(X, \shL_1^{d_1}\otimes\ldots\otimes\shL_r^{d_r})$.
Generalizing \emph{Grauert's Criterion} for ample sheaves,  we
characterize ample families also in terms of affine hulls and contractions.
Furthermore, we give a cohomological characterization which is analogous to
\emph{Serre's Criterion}
for ample sheaves. We also relate homogeneous spectra to
Cox's homogeneous coordinate rings for
toric varieties \cite{Cox 1995b}.

As an application, we shall generalize \emph{Grothendieck's Algebraization
Theorem}. The result is that a proper formal scheme
$\foX\ra\Spf(R)$ is algebraizable if
there is a finite collection of invertible formal sheaves restricting
to an ample family
on the closed fiber and satisfying an additional condition.

\begin{acknowledgement}
We thank Professor Uwe Storch for helpful suggestions.
The second author is grateful to the  M.I.T.\ Mathematical Department for
its hospitality, and thanks the DFG for financial support. We thank
Professor Paul Smith (University of Washington, Seattle) for pointing
out an error in the printed version and for showing us
Example \ref{greadedsurjectionrelevant}
\end{acknowledgement}

%===========================================================
\section{Grauert's criterion for ample families}

In this section, we shall generalize Grauert's criterion for ample sheaves to ample families.
Given a collection of invertible sheaves
$\shL_1,\ldots,\shL_r $,  we use
multiindices and set
$\shL^d=\shL_1^{d_1}\otimes\ldots\otimes\shL_r^{d_r}$  for   each
$d=(d_1,\ldots,d_r)\in\ZZ^r$. Let us start with the defining property of ample families:

\begin{proposition}
\label{ample families}
Let
$X$ be a quasicompact and
quasiseparated scheme. For a family
$\shL_1,\ldots,\shL_r$ of invertible sheaves, the following are equivalent:
\renewcommand{\labelenumi}{(\roman{enumi})}
\begin{enumerate}
\item The open sets
$X_f$ with
$f\in\Gamma(X,\shL^d)$ and
$d\in\NN^r$ form a base of the topology.
\item For each
$x\in X$, there is some $d\in\NN^r$ and
$f\in\Gamma(X,\shL^d)$ so that
$ X_f$ is an affine neighborhood of $x$.
\item For each point $x\in X$, there is a
$\QQ$-basis
$d_1,\ldots,d_r\in \NN^r$ and global sections
$f_i\in\Gamma(X,\shL^{d_i})$ so that
$ X_{f_i}$ are affine neighborhoods of $x$.
\end{enumerate}
\end{proposition}

\proof
To see
(ii) $\Rightarrow $ (iii), choose a degree
$d\in\NN^r$ and a section
$f\in\Gamma(X,\shL^{d})$ so that
$U=X_{f}$ is an affine neighborhood of the point
$x$. Choose a $\QQ$-basis $e_1,\ldots,e_r \in\NN^r$. Since
$U$ is affine, we find sections
$f_i' \in\Gamma(U,\shL^{e_i})$ satisfying
$f_i'(x)\neq 0$. According to \cite{EGA I} Theorem 6.8.1, the sections
$f_i=f_i'\otimes f^n$ extend from $U$ to
$X$ for
$n\gg 0$. The sections $f$ and $f_i$ have degrees
$d$ and $d_i=nd+e_i$, respectively. Skipping one of them we have a basis
and the desired sections.
The other implications are clear.
\qed

\medskip
Following Borelli \cite{Borelli 1963}, we call a finite collection
$\shL_1,\ldots,\shL_r$ of invertible
$\O_{X}$-modules an \emph{ample family}  if the scheme
$X$ is quasicompact and quasiseparated, and the equivalent conditions in Proposition \ref{ample families} hold.
A scheme  is called \emph{divisorial} if it admits an ample family of invertible sheaves.

Recall that a scheme is \emph{separated} if the diagonal embedding
$X\subset  X\times X$ is closed, and \emph{quasiseparated} if
the diagonal  is quasicompact.
Note that, in contrast to the definition of ample sheaves (\cite{EGA II} Def.\ 4.5.3), we do not require separateness
for ample families.
However, the possible nonseparatedness is rather mild:

\begin{proposition}
\label{affine diagonal}
The diagonal embedding of a divisorial scheme is  affine.
\end{proposition}

\proof
Let
$X$ be a divisorial scheme,
$\shL$ be an invertible sheaf,
$f\in\Gamma(X,\shL)$ a global section, and
$U=\Spec(A)$ an affine open subset. Then
$U\cap X_f=\Spec(A_f)$ is affine. Since
$X$ is covered by  affine open subsets of the form
$X_f$, this ensures that the diagonal embedding
$X\subset X\times X$ is affine.
\qed

\medskip
Obviously, schemes admitting  ample invertible sheaves are divisorial.
The following gives another large class of divisorial schemes:

\begin{proposition}
\label{Q-factorial schemes}
Normal noetherian locally
$\QQ$-factorial schemes  with affine diagonal   are divisorial schemes.
\end{proposition}

\proof
As in \cite{SGA 6} II 2.2.6, the complement of an affine dense open subset is a Weil divisor.
By assumption it is
$\QQ$-Cartier, so
$X$ is divisorial. Using quasicompactness, we find finitely many
effective Cartier divisors $D_1,\ldots,D_r\subset X$ with
$\bigcap D_i=\emptyset$. By Proposition \ref{ample families}, the
invertible sheaves $\O_X(D_i)$ form an ample family.
\qed

\medskip
Here are two useful properties of divisorial schemes:

\begin{proposition}
\label{resolution and affine torsor}
For divisorial noetherian schemes, the following hold:
\renewcommand{\labelenumi}{(\roman{enumi})}
\begin{enumerate}
\item Each coherent
$\O_{X}$-module admits a resolution with locally free
$\O_{X}$-modules of finite rank.
\item There is a noetherian ring
$A$ together with  a smooth surjective affine morphism
$\Spec(A)\ra X$.
\end{enumerate}
\end{proposition}

\proof
The first assertion is due to Borelli \cite{Borelli 1967} Theorem 3.3.
The second statement is called the  Jouanolou--Thomason trick. For a proof,
see \cite{Weibel 1989} Proposition 4.4.
\qed

\medskip
Grauert's Criterion
states that a line bundle is ample if and only if its zero section
contracts to a point  (see \cite{Grauert62} p.\ 341
and \cite{EGA II} Theorem 8.9.1).
The task now is to generalize this to families of line bundles.
To do so, we shall use vector bundles.
Recall that the category of locally free $\O_X$-modules $\shE$ is
antiequivalent to the category of vector bundles
$\pi:B\ra Z$ via $B=\Spec S(\shE)$ and
$\shE=\pi_*(\O_B)_1$. Under this correspondence,  the sections
$f\in\Gamma(X,S(\shE))$ correspond to functions $f\in\Gamma(B,\O_B)$.
For $f \in \Gamma(X,S^n(\shE))$ define
$X_f=\{ x \in X:\, f(x) \neq 0 \}$.
Then
$B_f \subset \pi^{-1}(X_f)$ and $X_f= \pi(B_f)$.

A locally free sheaf $\shE$ is called \emph{ample} if the invertible sheaf
$\O_P(1)$ is ample on $P=\PP(\shE)$ (see \cite{Hartshorne 1966}).
This easily implies that the open subsets $X_f\subset X$ with
$f\in\Gamma(X,S^n(\shE))$ generate the topology.
In contrast to line bundles, the latter condition is not sufficient for
ampleness. Let us characterize this condition:

\begin{theorem}
\label{Grauert criterion}
Suppose  $X$ is quasicompact and quasiseparated.
Let $\pi:B\ra X$ be a vector bundle, $\shE=\pi_*(\O_B)_1$
the corresponding locally free sheaf, and
$Z \subset B$ the zero section.
Then the following are equivalent.
\renewcommand{\labelenumi}{(\roman{enumi})}
\begin{enumerate}
\item
The open subsets
$X_f\subset X$ with $f\in\Gamma(X,S^n(\shE))$ generate the topology.
\item
For every point $x \in X$ there is a function $f \in \Gamma(B,\O_B)$
so that $B_f$ is affine and $\pi^{-1}(x) \cap B_f \neq \emptyset$.
\item
There is a scheme $B'$, and a morphism $q:B\ra B'$, and an open subset
$U' \subset B'$ so that the following holds:
The image $q(Z)\subset B'$ admits a quasiaffine open neighborhood,
$B \ra B'$ induces an isomorphism $q^{-1}(U') \cong U'$
and the projection $q^{-1}(U') \ra X$ is surjective.
\item
There exists an open subset $Z \subset W \subset B$ such that for every
$x \in X$ there is a function
$f \in \Gamma(W,\O_B)$ so that  $W_f$ is affine and
$W_f\cap \pi^{-1}(x)$ is nonempty.
\end{enumerate}
\end{theorem}
\proof
We shall prove the implications
(i) $\Rightarrow \ldots \Rightarrow$ (iv) $\Rightarrow$ (i).
First assume (i). Fix a point $x\in X$. Choose an affine neighborhood
$x \in V$ and a section
$f\in\Gamma(X,S^n(\shE))$ so that
$x \in X_f \subset V$.
Then $B_f =\pi^{-1}(V) _f$ is affine and $\pi^{-1}(x) \cap B_f \neq \emptyset$.

Assume (ii) holds. Set $B'=B^\aff = \Spec \, \Gamma(B,\O_B)$.
According to \cite{EGA I} Corollary 6.8.3,
the  map $\Gamma(B,\O_B)_f \ra \Gamma(B_f,\O_B)$
is bijective, so the affine hull $q: B \ra B^\aff$ induces an
isomorphism $B_f=q^{-1}(D(f)) \ra D(f)$ for $B_f$ affine.
Then take $U'=\bigcup_f D(f)$ and $U= \bigcup_f B_f$, where the union runs
over $f$, $B_f$ affine.

Assume that (iii) holds.
Let $ W' \subseteq B'$ be a quasiaffine neighborhood of $q(Z)$.
Then $W=q^{-1}(W')$ is an open neighborhood of the zero section.
Fix a point $x\in X$ and let $\eta\in\pi^{-1}(x)$ be the generic point,
such that $\eta \in q^{-1}(U') \cap W$.
Then we find $f\in \Gamma(W',\O_{B'})$ so that
$q(\eta) \in W'_f \subset U' \cap W'$ is affine
and $\eta \in W_f=q^{-1}(W'_f) \cong W'_f$ is an affine neighborhood.

Now suppose (iv) holds. Let $x \in V \subset X$. Then there exists
$f \in \Gamma(W,\O_B)$ so that $W_f \subset \pi^{-1}(V) \cap W $ is affine and
$\pi^{-1}(x) \cap W_f \neq \emptyset$.
Let $\shI\subset \O_W$ be the ideal of the zero section
$Z\subset W$.
Then $\O_W/\shI^{n+1}=\bigoplus_{d=0}^nS^d(\shE)$,
and $f \in \Gamma(W,\O_B)$ has a Taylor series expansion
$f=\sum_{d=0}^\infty f_d$ with $f_d \in \Gamma(X,S^d(\shE))$.
Choose a degree $d \geq 0$ with $f_d(x)\neq 0$.
Then $x \in X_{f_d} \subset V$
is the desired open neighborhood of $x$.
\qed

\begin{remark}
If $X$ is connected and proper over a base field, the image of the
zero section $Z\subset B$ in any quasiaffine scheme is a closed point.
In this case, the assumption in condition (iii)
implies that $q$ contracts $Z \cong X$ to a point.
\end{remark}

The preceding result yields a characterization of ample families
$\shL_1,\ldots,\shL_r$
in terms of the corresponding locally free sheaf
$\shE=\bigoplus\shL_i$.

\begin{corollary}
\label{ample family bundle}
Suppose
$X$ is   quasicompact and
quasiseparated. A family of invertible sheaves
$\shL_1,\ldots,\shL_r$ is ample if and only if the  vector
bundle $\pi:B\ra X$ with $\pi_*(\O_B)_1=\shL_1\oplus \ldots \oplus \shL_r$
satisfies the equivalent conditions
in Proposition \ref{Grauert criterion}.
\end{corollary}

\proof
If the family is ample, then condition (i) of Proposition
\ref{Grauert criterion} holds.
For the converse, fix a point $x\in X$, and choose a section
$f\in\Gamma(X,S^n(\shE))$ so that $x \in X_f \subset V$ where $V$ is affine.
Write $f=\sum f_d$ according to the decomposition
$S^n(\shE)=\bigoplus_d \shL^d$, where the sum
runs over all degrees $d\in\NN^r$ with $n=\sum d_i$. Pick a summand
$f_d$ with $f_d(x)\neq 0$.
Then $X_{f_d} \subset X_f \subset V$ is the desired affine open neighborhood.
\qed

\begin{remark}
In the situation of \ref{ample family bundle},
the vector bundle $B=B_1 \times_X \ldots \times_X B_r$
decomposes into $B_i = \Spec \, S(\shL_i)$. The affine hull
$B \ra B^\aff$ is an isomorphism outside
the coordinate hyperplanes.
The following examples illustrate what may happen on the union of the
coordinate hyperplanes.
\end{remark}

\begin{example}
Set $X=\PP^1\times\PP^1$. First, let $\shL_1$ and $\shL_2$ be invertible
sheaves  of bidegree $(1,0)$ and $(0,1)$,
respectively. This is an ample family because $\shL_1\otimes\shL_2$ is ample.
Set $B_i=\Spec \,S(\shL_i)$ and consider the
corresponding rank two vector bundle $B=B_1\times_X B_2$.
On each summand $B_i$,
the affine hull $B\ra B^\aff$ restrict to  the morphism
$$
B_i=\tilde{\AA}^2 \times\PP^1\stackrel{\pr_1}{\lra} \tilde{\AA}^2
\stackrel{g}{\lra} \AA^2
$$
where $g:\tilde{\AA}^2\ra \AA^2$ is the blowing-up of the origin.

Now let $\shL_1$ and $\shL_2$ be invertible sheaves of bidegree $(2,-1)$ and
$(-1,1)$,
respectively. This is an ample family, because $\shL_1^2\otimes\shL_2^3$
is an ample
sheaf. Here the affine hull $B\ra B^\aff$ contracts the union of the
coordinate hyperplane $B_1\cup B_2$ to a point.
\end{example}

In the next examples we consider the following situation.
Let $k$ denote a field and let $A$ be an $\NN$-graded $k$-algebra of
finite type with $A=k[A_1]$.
Let $n_1,\ldots,n_r \in \ZZ$ and consider the family
$\O_X(n_1),\ldots, \O_X(n_r)$ on $X= \Proj A$.
Now we  describe the corresponding
vector bundle and its ring of global sections.
(The second
description uses  $\Proj$ of a $\ZZ$-graded ring, which we introduce in the
next section.)

\begin{proposition}
Let $A, X$ and $n_1,\ldots ,n_r$ as above.
Let $S=A[T_1,\ldots,T_r]$ be $\ZZ$-graded by degrees $\deg(T_j)=-n_j$.
Then the vector bundle
$B=\Spec S(\O_X(n_1) \oplus \ldots \oplus \O_X(n_r))$ is
$$ B \cong \bigcup_{f \in A_1} D_+(f) \subset \Proj S_{\geq 0}
\quadand
B \cong \bigcup_{f \in A_1} D_+(f) \subset \Proj S \, .$$
If $A$ is normal and $\dim\, X \geq 1$, then $\Gamma(B,\O_B)= S_0$.
\end{proposition}
\proof
Let $R=S$ or $=S_{\geq 0}$.
$R_0$ has also a $\ZZ^r$-graduation with
$(R_0)_d = \{ aT^d: a \in A_{d_1n_1+\ldots+d_rn_r} \} $, $d=(d_1,\ldots,d_r)$.
There is a natural rational mapping $\Proj R \dasharrow \Proj A$ which is
defined on $\bigcup D_+(f)$, where the union runs over $f \in A_1$.
This is an affine morphism, thus we can check the identities
by looking at the rings of global sections and at the restriction maps.
We have
\begin{eqnarray*}
\Gamma(\pi^{-1}(D_+(f)), \O_B)_d
&=& \Gamma(D_+(f), S(\O_X(n_1) \oplus \ldots \oplus \O_X(n_r)))_{d} \cr
&=& \Gamma(D_+(f), \O_X(d_1n_1+\ldots +d_rn_r)) \cr
&=& (A_f)_{d_1n_1+\ldots +d_rn_r} \cr
&=& ((R_f)_0)_d \cr
&=& \Gamma(D_+(f), \O_{\Proj R})_d
\end{eqnarray*}
and the restriction maps respect these identities.
The last statement follows if we replace $D_+(f)$ by
$U=\bigcup_{f \in A_1} D_+(f)$. Then
$$\Gamma(U, \O_X(d_1n_1+\ldots +d_rn_r))
= \Gamma(D(A_1),(\O_{\Spec (A)})_{d_1n_1+\ldots +d_rn_r})
=A_{d_1n_1+\ldots +d_rn_r} \, ,$$
since $A$ is normal and $V(A_1)$ has codimension $\geq 2$.
\qed

\begin{example}
Let $X=\PP^m$ be the homogeneous spectrum of $k[X_0,\ldots, X_m]$, $m \geq 1$.
Let $\shL_1=\O_X(1)$ and $\shL_2=\O_X(-1)$.
Then $B_1=\Spec S(\O(1))$ is the blowing up of the vertex point
in $\AA^{m+1}$ and
$B_2=\Spec S(\O(-1))$ is the projection from a point in $\PP^{m+1}$.
The ring of global sections of the rank two vector bundle
$B=B_1 \times_X B_2$ is the polynomial algebra
$$
\Gamma(B,\O_B)=k[X_0S,\ldots,X_mS, ST] \subset k[X_0,\ldots,X_m][S,T]
$$
where $S,T$ are indeterminates with degrees $\deg(S)=-1$ and $\deg(T)=1$.
The affine hull $B\ra B^\aff$ contracts $B_2$
to a point and is an isomorphism
on the complement $B-B_2$. We have
$\Gamma(B|D_+(X_i),\O_B)= k[X_0/X_i, \ldots,X_m/X_i,X_iS,T/X_i]$,
and the affine hull is given by
$$
X_jS \mapsto \frac{X_j}{X_i}X_iS \quadand
ST\mapsto \frac{T}{X_i}X_iS .
$$
\end{example}
\begin{example}
Again let $A=k[X_0,\ldots,X_m]$ ($m \geq 1$) and consider
the family $\shL_1=\O_X(1), \ldots,\shL_r=\O_X(1)$.
Then $\Gamma(B,\O_B)$ is the determinantial algebra
$$
k[X_iT_j\mid 0\leq i\leq m, 1\leq j\leq r] \subset
k[X_0,\ldots,X_m,T_1,\ldots,T_r],
$$
generated by the entries of the $(m+1)\times r$-matrix $(X_iT_j)$,
with relations
given by the $2\times 2$-minors $(X_iT_j)(X_kT_l)-(X_iT_l)(X_kT_j)$.
The affine hull $\varphi:B\ra B^\aff$ contracts exactly the zero section
to a point.
The affine hull $\varphi$ may also be described as the blowing-up of the column ideal
$(X_0T_1,\ldots, X_mT_1)$.
For this blowing up is given by
$\Proj \, k[X_iT_j][X_0T_1U,\ldots,X_mT_1U]$, and
$$k[X_iT_j][X_0T_1U,...,X_mT_1U]
\cong k[X_iT_j][X_0,...,X_m]=k[X_0,...,X_m,T_1,...,T_r]_{\geq 0}\, .$$
\end{example}

%===========================================================
\section{The homogeneous spectrum of a multigraded ring}

Generalizing the classical notion of homogeneous coordinates,
Grothendieck defined   homogeneous spectra for
$\NN$-graded rings (\cite{EGA II} \S 2). In this section, we shall generalize his approach to multigraded rings.
Let
$D$ be a finitely generated abelian group and let
$S= \bigoplus_{d \in D} S_d$ be a
$D$-graded ring. Note that by \cite{SGA 3a} I  4.7.3, such gradings correspond to actions of the diagonizable
group scheme
$\Spec(S_0[D])$ on  the affine scheme
$\Spec(S)$.

According to geometric invariant theory (see \cite{Mumford; Fogarty; Kirwan 1993} Thm.\ 1.1), the projection
$\Spec(S)\ra\Spec(S_0)$ is a categorical quotient in the category of schemes.
There is a quotient
$\Spec(S)\ra\Quot(S)$ in the category of ringed spaces as well. In general,  the latter is quite
different from the first.
However, we have the following favorable situation. Call the ring
$S$  \emph{periodic} if the degrees of the homogeneous units
$f\in S^\times$ form a subgroup
$D'\subset D$ of finite index.
In this case, we may choose a free subgroup $D'\subset D$
of finite index,
such that  $S'= \bigoplus_{d \in D'} S_d$ is a
Laurent polynomial algebra
$S_0[T_1^{\pm 1},\ldots,T_r^{\pm 1}]$.

\begin{lemma}
\label{geometric quotient}
For periodic rings
$S$, the projection
$\Spec(S)\ra\Spec(S_0)$ is a geometric quotient in the sense of
geometric invariant theory.
\end{lemma}

\proof
Choose a free subgroup $D'\subset D$ of finite index as above. The corresponding inclusion of the
Veronese subring $S'\subset S$ is an integral ring extension, because $D/D'$ is torsion,  such that
$\Spec(S)\ra\Spec(S')$ is a closed morphism.
By \cite{Mumford; Fogarty; Kirwan 1993} Amplification 1.3, this morphism
is a geometric quotient.

Since $S'$ is a Laurent polynomial ring,
$\Spec(S')$ is a principal homogeneous space for the induced action of
$\Spec(S_0[\ZZ^r])$ and the projection
$\Spec(S')\ra\Spec(S_0)$ is a geometric quotient.
Being the composition of two geometric quotients,
$\Spec(S)\ra \Spec(S_0)$ is a geometric quotient as well.
\qed

\medskip
In light of this, we seek to pass from a given graded ring $S$ to periodic
rings via localization.
An element
$f\in S$ is called \emph{relevant} if it is homogeneous and  the localization
$S_f$ is periodic. Equivalently, the degrees of all homogeneous divisors
$g|f^n$, $n\geq 0$ generate a subgroup
$D'\subset D$ of finite index.
There exists then also an $n$ and a homogeneous factorization
$f^n=g_1 \cdots g_\ell$ such that the degrees $\deg(g_i)$ generate a
subgroup of finite index.
Since geometric quotients are quotients in the category of ringed spaces,
localization of relevant elements yields open subschemes
$$
D_+(f)=\Spec(S_{(f)})\subset \Quot(S)
$$
inside the ringed space
$\Quot(S)$. Here
$S_{(f)}\subset S_f$ is the degree zero part of the localization.
This leads to the following definition:

\begin{definition}
\label{homogeneous spectrum}
Let
$D$ be a finitely generated abelian group and let
$S=\bigoplus_{d\in D} S_d$ be a
$D$-graded ring. We define the scheme
$$
\Proj(S)= \bigcup_{f\in S  \text{ relevant}} D_+(f) \subset\Quot(S),
$$
and call it the \emph{homogeneous spectrum} of the
graded ring
$S$.
\end{definition}

For
$\NN$-gradings, this coincides with the usual definition.
As in the classical situation, we define
$S_+\subset S$ to be the ideal generated by all relevant
$f\in S$. The corresponding invariant closed subscheme
$V(S_+)$ is called the \emph{irrelevant subscheme}.
The complementary invariant open subset
$\Spec(S) -  V(S_+)$ is called the \emph{relevant locus}.
Obviously, we obtain an affine projection
$$\Spec(S) -  V(S_+)\lra \Proj(S),$$
which is a geometric quotient for the induced action.

\begin{remark}
The points  $x\in\Proj(S)$ correspond to  graded (not necessarily prime) ideals
${\fop} \subset S$ not containing
$S_+$ such that the subset of homogeneous elements
$H\subset S - \fop$ is closed under multiplication. The stalk of the structure sheaf at
$x\in\Proj(S)$ is canonically isomorphic to
$(H^{-1}S)_0$.
\end{remark}

To proceed, we need a finiteness condition for multigraded rings.
In the special case $D=\ZZ$, the following  is
due to Bruns and Herzog (\cite{Bruns; Herzog 1993} Thm.~1.5.5).

\begin{lemma}
\label{noetherian}
Let
$S$ be a ring graded by a finitely generated abelian group
$D$. Then the following are equivalent:
\renewcommand{\labelenumi}{(\roman{enumi})}
\begin{enumerate}
\item
The homogeneous ideals of
$S$ satisfy the ascending chain condition.
\item
The ring $S$ is noetherian.
\item
$S_0$ is noetherian and
$S$ is an
$S_0$-algebra of finite type.
\end{enumerate}
If $S$ is noetherian and $M \subset D$ is a finitely generated submonoid,
then $S_M$ is also noetherian.
\end{lemma}

\proof
The implications
(iii) $\Rightarrow $ (ii) and (ii) $\Rightarrow $(i) are trivial,
so assume that (i) holds.
We start with some preparations.
For a submonoid
$M\subset D$, let
$S_M=\bigoplus_{d\in M}S_d$ be the corresponding \emph{Veronese subring}.
Let $D' \subset D$ denote a subgroup.
Then the ring $S_{D'} \subset S$ is a direct summand. Therefore $S_{D'}$
satisfies also the ascending chain condition for homogeneous ideals.
In particular, $S_0$ is noetherian.
It follows at once that $S_e$ are noetherian $S_0$-modules.
Let $D \cong D' \oplus T$, where $T$ is finite and $D'$ free.
Then $S$ is finite over $S_{D'}$ and $S_{D'}$ fulfills (i).
Thus we may assume that $D$ is free.
Let us call a free submonoid
$M\subset D$ a
\emph{quadrant} if $M=\bigoplus\NN d_i$, where the $d_i\in D$
is a subset of a $\ZZ$-basis for $D$.

\begin{claim}
\label{finite}
For each quadrant $M\subset D$ and each degree $e\in D$, the $S_M$-module
$S_{M+e}=\bigoplus_{d\in M}S_{d+e}$ is finitely generated.
\end{claim}

We prove this by induction on $r=\rank{M}$.
Since $S_e$ are  noetherian $S_0$-modules, this holds for $r=0$.
Fix a quadrant $M$ of rank $r$ and
suppose the Claim is true for each quadrant of rank $r-1$.
Choose $d_1,\ldots,d_r\in D$ with $M=\bigoplus_{i=1}^{r}\NN d_i$, and let
$M_i=\bigoplus_{j\neq i}\NN d_j$ be the $i$-th boundary quadrant.

Condition (i) implies that the graded $S$-ideal $S\cdot S_{M+e}$ is finitely
generated. Choose  homogeneous generators $f_1,\ldots, f_s\in S_{M+e}$.
Then $S_d\subset \sum S_Mf_i$ for each $d\in M+e$ with
$d\geq \max\left\{\deg(f_1),\ldots, \deg(f_s)\right\}$.
Fix such a degree $d\in M+e$.
Clearly, there are finitely many $d_{ij}\in M+e$ with
$$
M+e = (M + d) + \sum_{i}\sum_j(M_i + d_{ij}).
$$
By induction, each $S_{M_i+d_{ij}}$ is a finitely generated $S_{M_i}$-module,
and
we conclude that $S_{M+e}$ is a finitely generated $S_M$-module.
This proves the claim.

\medskip
Fix a quadrant $M=\sum\NN d_i$ and set $M^*=M-0=\sum (M+d_i)$.
Then $S_{M +d_i}$ is a finitely generated $S_M$-module and
thus $S_{M^*} = \sum S_{M+d_i}$ is finitely generated.

Let $S_{M^*} =(f_1,\ldots,f_n)$, where $f_i$ are homogeneous
of degree $>0$. We show by induction on $M \cong \NN^r$ that
$S_M=S_0[f_1,\ldots,f_n]$.
Let $g \in S_d$, $d \in M$ and suppose that $S_e \subset S_0[f_1,\ldots,f_n]$
for all $e < d$.
We have $g = h_1f_1+ \ldots + h_nf_n$, where
$\deg (h_j)= \deg (g)- \deg (f_j) < d$ and the result follows.

Since $D$ is a finite union of quadrants with $\rank(M)=\rank(D)$,
we conclude
that $S$ is a finitely generated  $S_0$-algebra.

To proof the additional statement, let
$g_i, i \in I$ be a generating system for $S$
with degrees $d_i$.
We may assume that $M$ is saturated,
thus we may describe $M$ with finitely many linear forms $\psi_k: D \ra \ZZ$,
say $M = \bigcap _k \psi_k^{-1}(\NN)$.
Consider the mapping
$\varphi: \NN^I \ra D,\, (n_i)_{i \in I} \mapsto \sum_i n_id_i $.
Then $\varphi^{-1}(M)= \bigcap_k (\psi_k \circ \varphi)^{-1}(\NN) $
and therefore $\varphi^{-1}(M)$ is finitely generated,
say by $r_j, j \in J$. Set $e_j= \varphi (r_j)$ in $M$.
We claim that $S_M$ is generated by elements of degree $e_j, j \in J$.
An element $f \in S_d$, $d \in M$ can be written as
the sum of products $\prod_{i \in I} g_i^{n_i}$,
where $d= \sum _{i \in I} n_i d_i $.
But then $(n_i)_{i \in I} = \sum _j m_j r_j$ and
$$\prod_{i \in I} g_i^{n_i}
= \prod_{i \in I} g_i^{\sum_j m_j r_{ji} }
= \prod_{j \in J} ( \prod _{i \in I} g_i^{r_{ji}} )^{m_j} \, .$$
Thus $\prod_i g_i^{n_i}$ is a product of elements in $S_{e_j}$.
\qed

\medskip
We have the following finiteness condition for homogeneous spectra:

\begin{proposition}
\label{universally closed}
The morphism
$ \Proj(S) \ra\Spec(S_0)$ is universally closed and of finite type,
provided that
$S$ is noetherian.
\end{proposition}

\proof
By Lemma \ref{noetherian}, the ring $S_0$ is noetherian and the
$S_0$-algebra
$S$ is of finite type.
The relevant locus
$\Spec(S) -  V(S_+)$ is quasicompact and surjects onto
$\Proj(S)$, so the homogeneous spectrum is quasicompact.

Next, we check that the projection
$\Proj(S)\ra \Spec(S_0)$ is locally of finite type. To do so,
fix a relevant element
$f\in S$, and let
$D'\subset D$ be a free subgroup of finite index  such that for each
$d\in D'$, there is a homogeneous unit
$g\in S_f^\times$ with
$\deg(g)=d$. Let
$S'=\bigoplus_{d\in D'}S_d$ be the corresponding Veronese subring. By Lemma \ref{noetherian}, the ring
extension
$S'\subset S$ is finite, so Artin--Tate \cite{Artin; Tate 1951} tells us that the ring extension
$S_0\subset S'$ is of finite type.
Clearly, the localization
$S'_f$ is isomorphic to a Laurent polynomial ring
$S_{(f)}[T_1^{\pm 1},\ldots,T_r^{\pm 1}]$. Setting
$T_i=1$, we deduce that
$S_0\subset S_{(f)}$ is of finite type.

Finally, we verify universal closedness.
Let
$S_0 \ra R_0$ be a base change and set
$R=S\otimes_{S_0}R_0$.
By Lemma \ref{basechangelemma} below
we have
${\rm rad }(R_+) = {\rm rad}( S_+\otimes_{S_0}R_0)$ and so
$\Proj (S) \times_{\Spec (S_0)} \Spec (R_0)
= \Proj (S \otimes_{S_0}R_0)$.
Hence it is enough to check that
$h:\Proj(S)\ra \Spec(S_0)$ is closed under the new hypothesis that
$S$ is an
$S_0$-algebra of finite type, that
$S_+\subset S$ is an ideal of finite type, and that each
$S_d$ is a finitely generated
$S_0$-module. Each closed subset of $\Proj(S)$ is of the form
$\Proj( S/\fob)= V_+(\fob)$ for some graded ideal
$\fob\subset S$, and we have
$h(V_+({\fob})) \subset V({\fob}_0)$.
Consequently, it suffices to show that
$h:\Proj(S)\ra \Spec(S_0)$ has  closed image.

Fix a point
$x\in\Spec(S_0)$ with
$h^{-1}(x)=\emptyset$, and let
$\fop\subset S_0$ be the corresponding prime ideal. We  have to construct
$g\in S_0-\fop$ with
$h^{-1}(\Spec(S_0)_g)=\emptyset$. The condition
$h^{-1}(x)=\emptyset$ signifies that each relevant
$f\in S$ is nilpotent in
$ S/\fop S$. Hence
$S_+^k\subset \fop S_+$ for some integer
$k>0$. Choose finitely many homogeneous
$g_1,\ldots, g_n\in S$ with
$S=S_0[g_1,\ldots, g_n]$,  and set
$d_i=\deg(g_i)$. Call a degree
$d\in D$ \emph{generic} if for each linear combination
$d=\sum n_id_i$ with nonnegative coefficients, the set
$\left\{d_i\mid n_i\geq k\right\}$ generates a subgroup of finite index.
Then the set of nongeneric degrees
is a  union of a finite set and finitely many affine hyperplanes.

Let $d=\sum n_id_i$ be generic and consider
$g=\prod g_i^{n_i} \in S_d$. Then one may write
$g=g_1^k \cdot \ldots \cdot g_r^k \cdot g'$ such that
$g_1 \cdot \ldots \cdot g_r$ is relevant.
Using
$S_+^k\subset\fop S_+$, we see $g \in \fop S_d$.
This gives
$\fop S_d=S_d$ for all generic degrees
$d\in D$, and the Nakayama Lemma gives
$(S_d)_\fop=0$.

Next, choose finitely many relevant
$f_1,\ldots, f_m\in S$ with
${\rm rad}(S_+) \!=\! \sqrt{(f_1,\ldots, f_m)}$. We may assume
$f_i=f_{i1}\ldots f_{ir}$ such that each sequence
$\deg(f_{i1}),\ldots,\deg(f_{ir})\in D$ generates a subgroup of finite index.
For each
$f_i$, choose a linear combination with positive coefficients
$e_i=\sum_jn_j\deg(f_{ij})\in D$   that  is generic.
Then some
$g\in S_0-\fop$ annihilates all
$S_{e_i}$. Consequently, each
$f_i \in S[1/g]$ is nilpotent, hence the preimage
$h^{-1}(\Spec(S_0)_g)$ is empty.
\qed

\begin{remark}
Roberts (\cite{Roberts 1998}, sect.~8.2) introduced multihomogeneous
spectra for certain $\NN^r$-graded rings. To explain Roberts' conditions,
let $S_{[j]}\subset S$, $0\leq j\leq r$ be the graded subring generated by
all homogeneous elements whose degrees are of the
form $(i_1,\ldots,i_j,0,\ldots,0)$.
Then Roberts assumes that each $S_{[j+1]}$ is generated over $S_{[j]}$
by finitely many elements of degree $(i_1,\ldots,i_j,1,0,\ldots,0)$.
Now Roberts' homogeneous spectrum is the subset
$\bigcup D_+(f)\subset\Proj(S)$,
where the union runs over all $f\in S$ admitting a factorization
$f=g_1\ldots g_r$ so that $\deg(g_j)$ has the
form $(i_1,\ldots,i_{j-1},1,0,\ldots,0)$,
compare \cite{Roberts 1998}, Proposition 8.2.5.
\end{remark}
\medskip
The following lemma shows that the irrelevant locus behaves well
under base change.

\begin{lemma}
\label{basechangelemma}
Let
$S=\bigoplus_{d\in D} S_d$ be a
$D$-graded ring over
a finitely generated abelian group $D$ and let
$S_0 \rightarrow R_0$ be an arbitrary base change. Then we have
${\rm rad} (R_+) = {\rm rad} ((S_+) R)$,
where $R= S \otimes_{S_0}R_0$ and where $(S_+)R $ denotes the extended ideal under
the natural homomorphism $S \rightarrow R$.
\end{lemma}
\begin{proof} The inclusion
$\supseteq$ is clear, as a relevant element
$s \in S$ is also relevant after base change. So suppose that a homogeneous
element $r \in R$ of degree $d$ is relevant. Then we have a
factorization $r^k=r_1 \cdots r_\ell$ into homogeneous factors $r_i$
of degree $d_i$ such that the degrees
$d_i$ generate a subgroup of finite index.
For $i=1, \ldots , \ell$ we can write
$r_i= \sum_{j_i \in J_i}   s_{ij_i}\otimes a_{ij_i}$ with suitable index sets $J_i$, homogeneous elements
$s_{ij_i} \in S$ of degree $d_i$ and elements $a_{ij_i} \in R_0$.
Therefore
$$r^k
=(\sum_{j_1 \in J_1} s_{1j_1} \otimes a_{1j_1} ) \cdots
(\sum_{j_\ell \in J_\ell} s_{\ell j_\ell} \otimes a_{\ell j_\ell})
= \sum_{(j_1, \ldots , j_\ell) \in J_1 \times \ldots \times J_\ell}
 (\prod_{i=1}^\ell s_{ij_i})\otimes   (\prod_{i=1}^\ell a_{ij_i})\, .$$
For every choice $(j_1, \ldots ,j_\ell)$, the element $\prod_{i=1}^\ell s_{ij_i}$ is a product of elements of
degrees
$d_1, \ldots, d_\ell$ and so it is a relevant element in $S$. Therefore
$r^k \in (S_+) R$ and $r$ belongs to its radical.
\end{proof}

\begin{example}
\label{notrelevantbecomesrelevant}
The ideal generated by the relevant elements behaves not well under
base change, only its radical, and a homogeneous non-relevant
element might become relevant.
Consider
$S=K[u][x,y,z]/(yz-ux)$, where $K$ is a field, with the ${\mathbb Z}^2$-grading given by
$\deg (x)=e_1+ e_2$, $\deg (y)=e_1$, $\deg (z)=e_2$ and $S_0=K[u]$.
Consider the base change to the quotient field $K[u] \rightarrow K(u)=:R_0$.
The element $x$ is not relevant in $S$, since $S_x \cong K[x,y,z]_x$.
It becomes however relevant in
$R = S \otimes_{S_0} R_0$, since in this ring inverting $x$ makes also $y$ and $z$ to units.
The irrelevant ideals are given by
$S_+=(yz,xy,xz)$ and $R_+= (x)$, and $x=yz/u$.
\end{example}

%===========================================================
\section{Separation criteria}

In Proposition \ref{universally closed}, we may say that
$\Proj(S)$ is a \emph{complete}
$S_0$-scheme. We cannot, however, infer that
it is proper. For example, set
$S=k[X,Y]$ with degrees in
$D=\ZZ$ given by
$\deg(X)=1$ and
$\deg(Y)=-1$. Then
$S_0=k[XY]$ yields the affine line $\AA^1_k$, and
$\Proj(S)$ is the affine line with double origin, which is nonseparated.
However, the following holds:

\begin{proposition}
\label{affine diagonal homogeneous}
The diagonal embedding of a   homogeneous spectrum is affine.
\end{proposition}

 \proof
The intersection
$D_+(f)\cap D_+(g)=D_+(fg)$ is affine.
\qed

\medskip
Here is a criterion for separatedness, which trivially holds for
$\NN$-gradings:

\begin{proposition}
\label{separation criterion}
If for each pair
$x,y\in \Proj(S)$  there is a relevant
$f\in S$ with
$x,y\in D_+(f)$, then
$\Proj(S)$ is separated.
\end{proposition}

\proof
Under the assumption
the affine open subsets
$U\times U$, where
$U\subset \Proj(S)$ is affine, cover
$\Proj(S)\times\Proj(S)$. Clearly, the diagonal embedding of
$\Proj(S)$ is closed over each of these open subsets, hence it is closed.
\qed

\medskip
The next task is to recognize large separated open subsets in
$\Proj(S)$. Given a homogeneous element $f\in S$, let
$H_f\subset S$ be the set of homogeneous divisors
$g| f^n$,
$n\geq 0$, and
$C_f\subset D\otimes\RR$ the closed convex cone generated by  the degrees
$\deg(g)$, $g\in H_f$.
Note that a homogeneous element is relevant if and only if its cone  has nonempty interior.

\begin{proposition}
\label{separated subset}
Let
$f_i\in S$ be a collection of relevant elements so that  each closed convex cone
$C_{f_i}\cap C_{f_j}\subset D\otimes\RR$
has nonempty interior.
Then
$\bigcup D_+(f_i)\subset \Proj(S)$ is a separated open subset.
\end{proposition}

\proof
According to \cite{EGA I} Proposition 5.3.6, it suffices to show that the multiplication map
$S_{(f)}\otimes S_{(g)}\ra S_{(fg)}$ is surjective for each pair of relevant elements
$f,g\in S$ such that
$C_f\cap C_g$ has nonempty interior. Note that, for each factorization
$g^n=g_1\ldots g_m$, we may replace
$g$ by
$g_1^{n_1}\ldots g_m^{n_m}$,
$n_i>0$ without changing  the localization
$S_{(g)}$. Thus we may assume  $\deg (g)\in C_f$.
Passing to a suitable power of $g$, we may assume
$\deg(g)=\sum n_i\deg( f_i)$,
$n_i\geq 0$ with  $f_i \in H_f$.
Each element in
$S_{(fg)}$ has the form
$a/(fg)^k$ with
$a\in S$ homogeneous, so
$$
\frac{a}{(fg)^k} =
\frac{a}{f^k\prod f_i^{kn_i}}\cdot \frac{\prod f_i^{kn_i}}{g^k}
$$
is contained in the image of
$S_{(f)}\otimes S_{(g)}$.
\qed

\medskip
Next, we shall relate homogeneous spectra of multigraded polynomial
algebras to toric varieties. Fix a ground ring
$R$ and a free abelian group
$M$ of finite rank. A \emph{simplicial torus embedding} of the   torus
$T=\Spec(R[M])$ is an equivariant open embedding
$T\subset X$ that is locally given by
$R[M\cap\sigma^\vee]\subset R[M]$ for some strongly convex, simplicial cone
$\sigma\subset N_\RR$ in the dual lattice
$N=\Hom(M,\ZZ)$.
Here \emph{simplicial cone} means that the cone is generated by a linear independent set.
In contrast to the usual definition, we do not require that our torus embeddings are separated.

Simplicial torus embeddings occur in the following context,
which is related to a construction of Cox \cite{Cox 1995b}. Let
$S=R[T_1,\ldots, T_k]$ be a
$D$-graded polynomial algebra, such that the grading is given by a linear map
$d:\ZZ^k\ra D$ sending the
$i$-th base vector to
$\deg(T_i)\in D$. Let
$M\subset\ZZ^k$ be the kernel.

\begin{proposition}
\label{graded polynomial ring}
Notation as above. Then  $\Proj(S)$ is a  {\rm(}possibly nonseparated{\rm )}
simplicial torus embedding of
the torus
$\Spec(R[M])$.
\end{proposition}

\proof
Let
$I=\left\{ 1,\ldots ,k \right\}$ be the index set for the indeterminates.
Fix a relevant monomial
$T^n=T_1^{n_1}\ldots T_k^{n_k}$ and let
$J=\left\{ i\in I\mid n_i>0 \right\}$ be its support. A direct argument gives
$S_{(T^n)}=R[M_J]$ for the monoid
$$
M_J=(\ZZ^J\oplus\NN^{I-J})\cap M \subset \ZZ^k.
$$
Clearly, the submonoid
$M_J\subset M$ is the intersection of
$\Card(I-J)\leq\rank(M)$ half spaces. Therefore, its dual cone
$\sigma\subset N_\RR$ is simplicial.

It remains to check
$M=M_J+(-M_J)$.  Let
$D'\subset D$ be the subgroup generated by
$\deg(T_i)$ with
$i\in J$,  and
$m=\ord(D/D')$ be its index.
Then there are integers
$\lambda_j\in\ZZ$,
$j\in J$ solving the equation
$\sum_{i\in J}\lambda_i \deg(T_i)=-m\sum_{i\in I-J}\deg(T_i) $.
So the element
$g\in\ZZ^k$ defined by
$$
g_i=
\begin{cases}
\lambda_i & \text{for $i\in J$}\\
m & \text{for $i\in I-J$}
\end{cases}
$$
lies in
$M_J$. For each
$f\in M$, we have
$f+ng\in M_J$ for
$n\gg0$, hence
$f=(f+ng)-ng$  is contained in
$M_J+(-M_J)$.
\qed

\begin{corollary}
\label{divisorial}
If
$S$  is finitely generated as
$S_0$-algebra, then
$\Proj(S)$ is divisorial.
\end{corollary}

\proof
We may choose a
$D$-graded polynomial
$R$-algebra
$S'$ with
$S'_0=S_0$, together with a graded surjection
$S'\ra S$.
Such a graded surjection induces a closed embedding
$\Proj(S)\subset \Proj(S')$, because for every relevant element $f \in S$
we may find a relevant element $f' \in S'$ mapping to a power of it.
To see this let $f^n=g_1 \cdots g_\ell$ be a homogeneous factorization in $S$
such that the degrees of the $g_i$ generate a subgroup of finite
index. Let $g_1',..., g_\ell' \in S'$ be homogeneous elements
mapping to
$g_1,..., g_\ell$. Then their product is relevant in $S'$ and maps
to $f^n$.

By Proposition \ref{graded polynomial ring}, the scheme
$\Proj(S')$ is a simplicial torus embedding, which has affine diagonal by Proposition
\ref{affine diagonal homogeneous}, hence by Proposition
\ref{Q-factorial schemes} it is a
divisorial scheme. Consequently, the closed subscheme
$\Proj(S)$ is divisorial as well.
\qed

\begin{corollary}
\label{projective}
Suppose that $S$ is finitely generated over $S_0$.
If each finite subset of $\Proj(S)$ admits an affine neighborhood,
then $\Proj(S)$ is projective.
\end{corollary}

\proof
We already know that $\Proj(S)$ is of finite type, universally closed,
separated, and divisorial
(Prop.~\ref{universally closed}, Prop.~\ref{separation criterion},
and Cor.~\ref{divisorial}).
Since each finite subset admits an affine neighborhood,
the generalized Chevalley Conjecture (\cite{Kleiman 1966}, Thm.~3)
applies, and we conclude  that $\Proj(S)\ra\Spec(S_0)$ is projective.
\qed

\begin{remark}
Let us make the torus embedding
in Proposition \ref{graded polynomial ring} more explicit. For each subset
$J\subset I=\left\{ 1,\ldots,k \right\}$, let
$\sigma_J\subset N_\RR$ be the convex cone generated by the projections
$\pr_i:\ZZ^k\ra \ZZ$ restricted to $M$,
$i\in J$. You easily check that
$J\mapsto \sigma_{I-J}$
gives a bijection between the  subsets
$J\subset I$ with $\prod_{j \in J}T_j$ relevant,
and the strongly convex simplicial cones
$\sigma_{I-J}\subset N_\RR$. Let us call such subsets \emph{relevant}. Then the torus embedding is given by
$$
\Proj(S)=\bigcup_{J\subset I \text{ relevant}} \Spec(R[\sigma_{J-I}^\vee\cap M]).
$$
The (possibly nonseparated) union is taken with respect to the  inclusions
$J\subset J'$.
\end{remark}

\begin{example}
Let
$S=R[T_1,\ldots,T_k]$ be a polynomial algebra graded by
$D=\ZZ$ so that all indeterminates have positive degree. Then
$\Proj(S)$ is   the weighted projective space studied by
Delorme \cite{Delorme 1975}, Mori \cite{Mori 1975}, and
Dolgachev \cite{Dolgachev 1981}.
\end{example}

\begin{example}
Here we construct a separated non-quasiprojective scheme
defined by a single equation inside a multihomogeneous spectrum.
Let
$$
S=k[X_1,\ldots, X_4,Y_1,\ldots,Y_4,Z]
$$
be a $\ZZ^2$-graded polynomial ring over a field $k$,
with degrees
$\deg(X_i)=(1,0)$, $\deg(Y_j)=(0,1)$, and $\deg(Z)=(1,1)$.
Set $P=\Proj(S)$ and consider the open subset
$U=D_+(X_1Z)\cup D_+(Y_1Z)$. This is not separated:
We have
\begin{equation}
\label{local sections}
\Gamma(D_+(X_1Z),\O_P) = k[\frac{X_i}{X_1},\frac{X_1Y_m}{Z}]\quadand
\Gamma(D_+(Y_1Z),\O_P) = k[\frac{Y_j}{Y_1},\frac{X_mY_1}{Z}],
\end{equation}
with $1\leq i,j,m\leq 4$.
On the intersection $D_+(X_1Y_1Z)=D_+(X_1Z)\cap D_+(Y_1Z)$,
these algebras generate the subalgebra
$$
k[\frac{X_i}{X_1},\frac{Y_j}{Y_1},\frac{X_1Y_1}{Z}]\subset
k[\frac{X_i}{X_1},\frac{Y_j}{Y_1},\frac{X_1Y_1}{Z},\frac{Z}{X_1Y_1}]
=\Gamma(D_+(X_1Y_1Z),\O_P),
$$
which does not contain $Z/X_1Y_1$.
To obtain separated subschemes, we have to kill $Z/X_1Y_1$.
Consider the homogeneous polynomials of degree $(2,2)$
$$
g= X_1Y_1Z + X_2^2Y_1^2 + X_1^2Y_2^2 \quadand
f= X_1Y_1Z + X_2^2Y_1^2 + X_1^2Y_2^2 + X_3X_4Y_3Y_4.
$$
Let $S'=S/(f)$,
$P'=\Proj(S')$ and $U'= U\cap P'$.
Modulo $f$, the element
$Z/X_1Y_1$
is generated by the algebras in (\ref{local sections}), thus
$U'$ is a separated scheme. It is, however, not quasiprojective.
First observe that $S'$ is a factorial domain:
$S'_{X_3X_4Y_3} \cong k[X_0,\ldots,X_4,Y_1,Y_2,Y_3,Z]_{X_3X_4Y_3}$
is factorial and $X_3,X_4,Y_3$ are prime in $S'$, because
$g \in k[X_1,X_2,Y_1,Y_2,Z]$ is prime.

Choose points $x\in V_+(X_1,\ldots,X_4)\cap U'$
and $y \in V_+(Y_1,\ldots,Y_4)\cap U'$ (such points exist), and assume that
they admit a common affine neighborhood $W\subset U'$.
Then the preimage $V\subset\Spec(S')$ is affine as well.
By factoriality, $V=D(h)$ for some
homogeneous $h\in S$ with $h\in (X_1Z,Y_1Z)$.
Write $h=pX_1Z+qY_1Z$.
Let $\deg (h) = (d_1,d_2)$ and suppose $d_1 \geq d_2$.
Since $Y_1Z$ has degree $(1,2)$, it follows that
$q \in (X_1,\ldots,X_4)$. But then $h(x)=0$, contradiction.

\end{example}

\begin{example}
\label{greadedsurjectionrelevant}
For a graded surjection $T \rightarrow S$, a relevant
element $ s \in S$ does not necessarily come from a relevant element
in $T$ (only some power of it, as shown in
Corollary \ref{divisorial}).
The following example was communicated to us by Paul Smith.
Let $T=K[x,y,z]$ be ${\mathbb Z}^2$-graded with degrees $\deg(x)=e_1$,
$\deg(y)=-e_1$ and $\deg (z)=e_2$. Consider $S:=T/((xy-1)z^2)$ with
the induced grading. Then $z$ is relevant in $S$, since $z^2=xyz^2$.
However, $z$ is not the image of a relevant element in $T$.
Such an element $g$ must have degree $e_2$ and be of the form
$g=z+(xy-1)z^2h$ for some $h \in T$.
Since $\deg(xy-1)z^2=2e_2$, and since there are no elements of degree $-e_2$, it follows that $h=0$.
So $g=z$, but the homogeneous units in $K[x,y,z,z^{-1}]$ have all their degrees in ${\mathbb Z} e_2$.
\end{example}

%===========================================================
\section{Ample families and mappings to homogeneous spectra}

In this section, we shall relate homogeneous spectra to ample families.
Let
$X$ be a scheme,
$D$ a finitely generated abelian group, and
$\shB=\bigoplus_{d\in D}\shB_d$ a quasicoherent
$D$-graded
$\O_{X}$-algebra. We say that
$\shB$ is \emph{periodic} if each stalk
$\shB_x$ is a  periodic
$\O_{X,x}$-algebra. For each homogeneous
$f\in\Gamma(X,\shB)$, let
$X_f\subset X$ be the largest open subset such that all   multiplication maps
$f^n:\O_{X}\ra\shB_{nd}$ with
$n\geq 0, d=\deg(f)$ are bijective.

\begin{proposition}
\label{periodic ring}
Let
$S$ be a
$D$-graded ring,
$X=\Proj(S)$ its homogeneous spectrum,
$Y=\Spec(S) -  V(S_+)$ the relevant locus, and
$\pi:Y \ra X$ the natural projection. Then
$\pi_*(\O_{Y})$ is a periodic
$\O_X$-algebra. Furthermore, $X_f=D_+(f)$ for each relevant
$f\in S$.
\end{proposition}

\proof
By definition,
$S_g$ is a periodic
$S_{(g)}$-algebra for each relevant
$g\in S$, so
$\pi_*(\O_{Y})$ is a periodic
$\O_{X}$-algebra with
$\pi_*(\O_{Y})_0=\O_{X}$.
The inclusion
$D_+(f)\subset X_f$ is obvious. To verify
$X_f\subset D_+(f)$, it suffices to check that
$f\in S_g$ is invertible for each relevant
$g\in S$ with
$D_+(g)\subset X_f$. Replacing
$f$ by  a positive multiple and
$S_g$ by a suitable Veronese subring, the ring
$S_g$ becomes isomorphic to the Laurent polynomial algebra
$S_{(g)}[T_1^{\pm 1},\ldots,T_r^{\pm 1}]$, and
$f$ corresponds to a monomial
$\lambda T^d$ with
$\lambda\in S_{(g)}^\times$. Hence
$f$ is invertible.
\qed

\medskip
Next, we extend Grothendieck's description (\cite{EGA II} Prop.~3.7.3)
of mappings into homogeneous spectra to the multigraded case:

\begin{proposition}
\label{maps to homogeneous spectra}
Let
$X$ be a scheme,
$\shB$ a quasicoherent
$D$-graded
$\O_{X}$-algebra,
$S$ a
$D$-graded ring, and
$\varphi:S\ra  \Gamma(X,\shB)$ a graded homomorphism. Set
$U=\bigcup X_{\varphi(f)}$, where the union runs over all relevant
$f\in S$. Then there is a natural
morphism
$r_{\shB,\varphi}:U\ra \Proj(S)$ and a commutative diagram
$$
\begin{CD}
U @<<< \Spec(\shB_U)     @>>> \Spec(\shB)    \\
@Vr_{\shB,\varphi} VV @VV\varphi V     @VV\varphi V\\
\Proj(S)  @<<< \Spec(S) -  V(S_+)  @>>>  \Spec(S).
\end{CD}
$$
\end{proposition}

\proof
Each relevant
$f\in S$ gives a homomorphism
$S_{(f)}\ra \Gamma(X,\shB)_{(f)}$, where we write  $f$ instead of $\varphi(f)$.
Furthermore, there is a homomorphism
$$
\Gamma(X,\shB)_{(f)}\lra \Gamma(X_f,\O_X),\quad g/f^n\mapsto (f^n|X_f)^{-1}(g),
$$
where
$(f^n|X_f)^{-1}$ is the inverse mapping for the bijective  multiplication mapping
$f^n|X_f:\O_{X}|X_f\ra\shB_{nd}|X_f$. The composition defines a morphism
$X_f\ra D_+(f)$. You easily check that these morphisms coincide on the overlaps,
and we obtain the desired morphism
$r_{\shB,\varphi}:U\ra \Proj(S)$.
\qed

\medskip
We write
$r_{\shB,\varphi}:X\dasharrow\Proj(S)$ for the morphism
$r_{\shB,\varphi}:U\ra \Proj(S)$ and call it a \emph{rational map}. Saying that  a rational map is
\emph{everywhere defined} means
$U=X$. In this case, we have a honest morphism
$r_{\shB,\varphi}:X\ra \Proj(S)$.

\begin{corollary}
\label{periodic algebra via map}
Let
$S$ be a
$D$-graded ring. For each morphism
$r:X\ra \Proj(S)$, there is a quasicoherent periodic $D$-graded
$\O_{X}$-algebra
$\shB$ and a homomorphism
$\varphi:S\ra\Gamma(X,\shB)$ such that the rational map
$r_{\shB,\varphi}:X\dasharrow \Proj(S)$ is everywhere defined and coincides with
$r:X\ra \Proj(S)$.
\end{corollary}

\proof
Let
$Y=\Spec(S) -  V(S_+)$ be the irrelevant locus,
$\pi:Y\ra \Proj(S)$ the canonical projection, and set
$\shB=r^*(\pi_*(\O_{X}))$.
\qed

\medskip
We come to the characterization of ample families in terms of homogeneous
spectra:

\begin{theorem}
\label{families and spectra}
Let
$\shL_1,\ldots, \shL_r$ be a family of invertible
sheaves on a quasicompact and quasiseparated scheme
$X$. Then the following conditions are equivalent:
\renewcommand{\labelenumi}{(\roman{enumi})}
\begin{enumerate}
\item The family
$\shL_1,\ldots, \shL_r$ is ample.
\item The canonical rational map
$X\dasharrow \Proj \Gamma(X,\bigoplus_{d\in\NN^r}\shL^d)$ is everywhere defined and an open embedding.
\end{enumerate}
If
$X$ is of finite type over a noetherian ring
$R$, this is also equivalent with
\renewcommand{\labelenumi}{(\roman{enumi})}
\begin{enumerate}
\item[(iii)] There is a finite family of sections
$f_i \in \Gamma(X,\shL^{d_i}), i\in I$ and a
$D$-graded polynomial algebra $A=R[T_i]_{ i \in I}$
such that the rational map
$X\dasharrow \Proj(A)$ induced by $T_i \mapsto f_i$ is everywhere defined and an embedding.
\end{enumerate}
\end{theorem}

\proof
Set
$S=\Gamma(X,\bigoplus_{d\in\NN^r}\shL^d)$. We start with the implication
(i) $\Rightarrow$ (ii).
According to Lemma \ref{ample families}, for each point
$x\in X$, there is a
$\QQ$-basis
$d_i\in\NN^r$ and sections
$f_i\in \Gamma(X,\shL^{d_i})$ so that
$X_{f_i}$ are affine neighborhoods of
$x$. Consequently,
$f=f_1\ldots f_r\in S$ is relevant, so the rational map
$X\dasharrow \Proj(S)$ is everywhere defined.
Fix a relevant
$f\in S$ so that
$X_f$ is  affine.  According to \cite{EGA I}, the canonical map
$\Gamma(X,\bigoplus\shL^d)_f\ra \Gamma(X_f,\bigoplus\shL^d)$ is bijective.
Consequently,
$X_f\ra D_+(f)$ is an isomorphism, so
$X\ra \Proj(S)$ is an open embedding.
The reverse implication is trivial.

For the rest of the proof, suppose that
$X$ is of finite type over a noetherian ring
$R$. Assume that
$\shL_1,\ldots, \shL_r$ is ample. Choose finitely many relevant
$f_i\in S$ so that
$X_{f_i}\subset X$ form an affine open cover.
Due to the assumption we may write
$\Gamma(X_{f_i},\O_{X})=R[h_{i 1},\ldots,h_{i m_i}]$.
For suitable $n$ we have $f_{ij}:=f^{n}_i h_{ij} \in S $.
Let $f_i \in S, i \in I$ be these elements all together.

Let $R[T_i]$ be graded by $d(T_i)=\deg (f_i)$ such that the natural mapping
$R[T_i] \ra S$ is homogeneous.
Then
$X_{f_i} \ra D_+(f_i)$ are closed embeddings, since the ring morphisms are
surjective, and so
$X \ra \Proj(A)$ is an embedding.
Finally, the implication
(iii) $\Rightarrow$ (ii) is trivial.
\qed

\begin{corollary}
\label{quasiaffine}
Let
$\shL_1,\ldots,\shL_r$ an ample family on $X$. For each relevant element
$f\in \Gamma(X,\bigoplus_{d\in\NN^r}\shL^d)$, the open subset
$X_f$ is quasiaffine.
\end{corollary}

\proof
According to Theorem \ref{families and spectra}, we have an open embedding
$X_f\subset D_+(f)$.
\qed

\begin{corollary}
\label{isomorphism}
Let
$\shL_1,\ldots,\shL_r$ an ample family of invertible sheaves on $X$.  Set
$S=\Gamma(X,\bigoplus_{d\in\NN^r}\shL^d)$.
Then the following conditions are equivalent:
\renewcommand{\labelenumi}{(\roman{enumi})}
\begin{enumerate}
\item The open embedding
$r:X\ra \Proj(S)$ is an isomorphism.
\item For each relevant
$f\in S$, the quasiaffine open subset
$X_f\subset X$ is affine.
\end{enumerate}
If the affine hull
$X\ra X^\aff$ is proper, this is also equivalent to
\renewcommand{\labelenumi}{(\roman{enumi})}
\begin{enumerate}
\item[(iii)] The homogeneous spectrum
$\Proj(S)$ is separated.
\end{enumerate}
\end{corollary}

\proof The equivalence
(i) $\Leftrightarrow$ (ii) follows from $D_+(f)=X_f^\aff$.
Now  assume that
$X\ra X^\aff$ is proper.
The implication
(i) $\Rightarrow$ (iii) is trivial.
To see the converse, we apply \cite{EGA II} Corollary 5.4.3  and infer that the open dense
embedding
$X_f\ra D_+(f)$ is  proper, hence an isomorphism.
\qed

\medskip
Finally, we generalize Hausen's \cite{Hausen 2000} characterization of    divisorial varieties:

\begin{corollary}
\label{characterization of divisorial}
Let
$X$ be a scheme of finite type over a noetherian ring
$R$. Then the following are equivalent:
\renewcommand{\labelenumi}{(\roman{enumi})}
\begin{enumerate}
\item The scheme
$X$ is divisorial.
\item
There is an embedding of
$X$ into the homogeneous spectrum of a
multigraded
$R$-algebra of finite type.
\item
$X$ is embeddable into a simplicial torus embedding with affine diagonal.
\end{enumerate}
\end{corollary}

\proof
If
$X$ is divisorial,  Theorem \ref{families and spectra} ensures the existence of
an embedding
$X\subset \Proj(S)$ with
$S$ finitely generated.  The implication
(ii) $\Rightarrow$ (iii) follows from Proposition
\ref{graded polynomial ring}, and
(iii) $\Rightarrow$ (i) is trivial.
\qed

%===========================================================
\section{Cohomological characterization of ample families}

Throughout this section,
$R$ is a noetherian ring and
$X$ is a proper
$R$-scheme. According to Serre's Criterion (\cite{EGA IIIa} Prop.~2.6.1),
an invertible
$\O_{X}$-module
$\shL$ is ample if and only if for each coherent
$\O_{X}$-module
$\shF$ there is an integer
$n_0$ so that
$H^p(X,\shF\otimes \shL^n)=0$ for all
$p>0$, $n>n_0$. The task now is to generalize this to ample families.

Let $\shA=\bigoplus_{n\geq 0}\shA_n$ be a graded quasicoherent
$\O_X$-algebra of finite type generated
by $\shA_1$. Then $P=\Proj(\shA)$ is a projective $X$-scheme and $\O_P(1)$
is an
$X$-ample invertible sheaf. In general, however, $\O_P(1)$ is not ample
in the absolute sense. More precisely:

\begin{proposition}
\label{ample blowing-up}
With the preceding notation, the invertible
sheaf
$\O_P(1)$ is ample if and only if for each coherent
$\O_{X}$-module
$\shF$, there is an integer
$n_0$ so that
$H^p(X,\shF\otimes\shA_n)=0$ for
$p>0$ and  $n>n_0$.
\end{proposition}

\proof
Let
$h:P\ra X$ be the canonical projection. First, suppose that
$\shM=\O_P(1)$ is ample. Choose
$n_0>0$ so that  the canonical map
$\shA_n\ra h_*(\shM^n)$ is bijective and that
$R^qh_*(h^*(\shF)\otimes\shM^n)=0$ and
$H^p(P,h^*(\shF)\otimes\shM^n)=0$ holds for
$p,q>0$,
$n>n_0$. Using the Leray--Serre spectral sequence we infer
$H^p(X,h_*(h^*(\shF)\otimes\shM^n))=0$ for all
$p>0$,
$n>n_0$.

We claim that  the adjunction map
$\shF\otimes h_*(\shM^n)\ra h_*(h^*(\shF)\otimes\shM^n)$ is bijective for
$n\gg0$.  Fix a point
$x\in X$ and choose a finite presentation
$$
\O_{X,x}^{\oplus k}  \lra \O_{X,x}^{\oplus l}   \lra \shF_x  \lra 0.
$$
Then we have  an exact sequence
$$
\O_{X,x}^{\oplus k}\otimes h_*(\shM^n) \lra\O_{X,x}^{\oplus l}\otimes h_*(\shM^n) \lra\shF_x\otimes h_*(\shM^n)
\lra0,
$$
and another exact sequence
$$
h_*(\O_{P}^{\oplus k}\otimes \shM^n)    \lra  h_*(\O_{P}^{\oplus l}\otimes \shM^n) \lra h_*(h^*(\shF)\otimes\shM^n)
\lra  R^1h_*(\shG\otimes\shM^n)
$$
on $\Spec(\O_{X,x})$, where
$\shG$ is the kernel of
$\O_{P}^{\oplus l}\ra h^*(\shF)$. But
$R^1h_*(\shG\otimes\shM^n)_x=0$ for
$n\gg0$. By the 5-Lemma, the Claim is true locally around
$x$. Using  quasicompactness, we infer that the Claim holds  globally.
Enlarging
$n_0$ if necessary, we have
$H^p(X,\shF\otimes\shA_n)=0$ for
$p>0$ and
$n>n_0$ as desired. The converse is similar.
\qed

\medskip
Let $\shL$ be an invertible $\O_X$-module.
The idea now is to consider coherent submodules
$\shK\subset \shL$. Note that such submodules correspond to quasicoherent
graded subalgebras
$\bigoplus_{n\geq 0}\shK^n\subset S(\shL)$ locally
generated by terms of degree one.
In turn, we obtain a projective
$X$-scheme
$P=\Proj(\bigoplus_{n\geq 0}\shK^n)$ endowed with an
$X$-ample invertible sheaf
$\O_P(1)$. Locally,  $P\ra X$ looks like  the blowing-up of an ideal.

\begin{proposition}
\label{local topology}
Let
$\shL$ be an invertible
$\O_{X}$-module, and
$x\in X$ a closed point.  Then the following are equivalent:
\renewcommand{\labelenumi}{(\roman{enumi})}
\begin{enumerate}
\item For some
$n>0$, there is  a section
$f\in H^0(X,\shL^n)$ so that
$  X_f$ is an affine neighborhood of
$x$.
\item
For some
$d>0$, there is a coherent submodule
$\shK\subset\shL^d$ with
$\shK_x\subset \shL^d_x$ bijective, so that the graded $\O_X$-algebra
$\shA=\bigoplus_{n\geq 0}\shK^n$ satisfies the equivalent
conditions of Proposition \ref{ample blowing-up}.
\end{enumerate}
\end{proposition}

\proof
First, we check
(i) $\Rightarrow$ (ii). Set
$U=X_f$. According to \cite{Luetkebohmert 1990} Proposition 5.4, there is a blowing-up
$h:P\ra X$ with center disjoint from
$U$, together with an effective ample  Cartier divisor
$D\subset P$ satisfying
$\Supp(D)=P -  U$. Set
$\shM=\O_P(D)$, such that
$P=\Proj(\bigoplus_{n\geq 0} h_*(\shM^n))$. Replacing
$D$ by a suitable multiple, we may assume that
$\bigoplus_{n\geq 0} h_*(\shM^n)$ is generated by terms of degree one.

The identity
$g_U : h_*(\shM)|U \ra \O_{X}|U$ extends to a mapping
$g: h_*(\shM) \ra \shL^d $ for
$d\gg 0$.
This map is injective because a nonzero  section
$t \in \Gamma(P,\O_P(D))=\Gamma(X,h_*(\shM))$
does not vanish on $U$.
Let $\shK \subset \shL^d$ be the image of $g: h_*(\shM) \ra \shL^d$. Then
$\shK^n=h_*(\shM^n)$, because
$\bigoplus_{n\geq 0} h_*(\shM^n)$ is generated by terms of degree one. On
$P=\Proj(\bigoplus \shK^n)$ we have
$\O_P(1)= \shM$, which is ample as desired.

Now we check (ii) $\Rightarrow$ (i). Set $P=\Proj(\shA)$ and
let $h:P\ra X$ be the canonical morphism. Choose an affine open
neighborhood $U\subset X$ of $x$ so that the induced
projection $h^{-1}(U)\ra U$ is an isomorphism.
Since $\O_P(1)$ is ample, there is an integer $m>0$
and a section $g\in H^0(X,\shK^m)$ so that the induced section
$g'\in H^0(P,\O_P(m))$
vanishes on $P-h^{-1}(U)$ and is nonzero on the point  $h^{-1}(x)$.
Let $f$ be the image of $g$ under the inclusion
$H^0(X,\shK^m) \hookrightarrow H^0(X,\shL^{dm})$.
Let $y \in X_f$. Then $\shK^m_y = \shL^{dm}_y$ and $X_f$ lies inside
the locus where $h: P \ra X$ is an isomorphism.
Therefore
$P_g \cong X_f \subset U$ and $X_f$ is affine.
\qed

\medskip
The preceding result leads to a characterization of ample families in terms of cohomology:

\begin{theorem}
\label{families and cohomology}
Let
$R$ be a noetherian  ring and
$X$ a proper
$R$-scheme. A family
$\shL_1,\ldots,\shL_r$ of invertible
$\O_{X}$-modules is ample if and only if for each closed point
$x\in X$, there is
$d\in\NN^r$ and a coherent subsheaf
$\shK\subset \shL^d$ with
$\shK_x=\shL_x^d$ so that for each coherent
$\O_{X}$-module
$\shF$ there is an integer
$n_0$ with
$H^p(X,\shF\otimes\shK^n)=0$ for all
$p>0$, $n>n_0$.
\end{theorem}

\proof
This follows directly from  Proposition \ref{local topology}.
\qed

%===========================================================
\section{Algebraization of formal schemes via ample families}

Throughout this section,
$(R,\fom,k)$ denotes a complete local noetherian  ring, and
$\foX\ra \Spf(R)$ is a proper formal scheme.
Such formal schemes frequently occur as formal solutions for problems related
to moduli spaces and deformation theories. A natural question is whether such a formal scheme is \emph{algebraizable}.
This means that there is a proper scheme
$X\ra \Spec(R)$ whose
$\fom$-adique completion  is isomorphic to
$\foX$.  Grothendieck's Algebraization Theorem (\cite{EGA IIIa} Thm.\  5.4.5) asserts that
$\foX$ is algebraizable if there is an invertible
$\O_{\foX}$-module whose restriction to the closed  fiber
$X_0=\foX\otimes k$ is ample.
Here is a  generalization:

\begin{theorem}
\label{algebraization}
Let $\foX$ be a proper formal scheme as above, and
$\shL_1,\ldots,\shL_r$ a family of invertible
$\O_{\foX}$-modules.
Suppose that for each closed point $x\in\foX$, there is a degree
$d\in\NN^r$ and a coherent submodule $\shK\subset\shL^d$
with $\shK_x\subset\shL^d_x$ bijective, so that $\O_{P_0}(1)$ is ample
on $P_0=\Proj(\bigoplus_{m\geq 0}\shK^m/\shI\shK^m)$. Then  the formal scheme
$\foX$ is algebraizable.
\end{theorem}

\proof
Set
$X_n=\foX\otimes R/\fom^{n+1}$, such that
$\foX=\limdir X_n$, and let $\shI\subset\O_{\foX}$ be the
ideal of the closed fiber $X_0\subset\foX$.
As in the proof of Proposition \ref{local topology}, there is an integer $m>0$
and a global section $s_0\in H^0(X_0,\shK^m/\shI\shK^m)$
so that the induced section $t_0\in\Gamma(X_0,\shL^{dm}/\shI\shL^{dm})$
defines an affine open neighborhood $x\in (X_0)_{t_0}$.
We seek to extend such sections to formal sections.

For each $m\geq 0$, set
$\shA_m=\bigoplus_{n\geq 0} \shI^n\shK^m/\shI^{n+1}\shK^m$,
and let $\shA=\bigoplus_{m\geq 0}\shA_m$ be the corresponding $\NN$-graded
quasicoherent $\O_{X_0}$-algebra. Consider its  homogeneous spectrum
$P=\Proj(\shA)$.
We claim that the invertible sheaf  $\O_{P}(1)$ is  ample.
To see this, let $X'$ be the affine $X_0$-scheme defined by
$\O_{X'}=\bigoplus_{n\geq 0}\shI^n/\shI^{n+1}$.
Note that $X'$ is proper over the noetherian ring
$R'=\bigoplus\fom^n/\fom^{n+1}$.
Set
$$
P'=\Proj(\bigoplus_{m\geq 0}
(\bigoplus_{n\geq 0} \shI^n/\shI^{n+1} \otimes \shK^m/\shI\shK^m )),
$$
where the homogeneous spectrum is taken with respect to the grading $m\in\ZZ$.
We have $P'=X' \times_{X_0}P_0$ and conclude that $\O_{P'}(1)$ is ample.

The surjective mapping
$\shI^n \otimes_\foX \shK^m \ra \shI^n\shK^m$ induces a surjective
mapping
$$\shI^n/\shI^{n+1} \otimes_{X_0} \shK^m/ \shI \shK^m =
\shI^n \otimes_\foX \shK^m \otimes_\foX \O_{X_0}
\lra \shI^n\shK^m \otimes_\foX \O_{X_0} =\shI^n\shK^m/\shI^{n+1}\shK^m \, .$$
This yields a closed embedding $P\subset P'$,
showing that $\O_P(1)$ is ample.
To proceed, consider the exact sequence
$$
H^0(X_{n},\shK^m/\shI^{n+1}\shK^m ) \lra H^0(X_{n-1},\shK^m/\shI^{n}\shK^m )
\lra H^1(X_0,\shI^n\shK^m/\shI^{n+1}\shK^m).
$$
By Proposition \ref{ample blowing-up}, there is an integer $m_0>0$ so that
the group on the right is zero for all $m\geq m_0$ and all $n\geq 0$.
Passing to a suitable multiple if necessary, we can lift
our section $s_0\in H^0(X_0,\shK^m/\shI\shK^m)$ to a formal section
$s\in H^0(\foX,\shK^m)$.
Therefore, the section
$t_0\in\Gamma(X_0,\shL^{dm}/\shI\shL^{dm})$ lifts to a formal section
$t\in\Gamma(\foX,\shL^{dm})$.

Using such formal sections, you construct as in the proof of Theorem
\ref{families and spectra}  a finitely generated
$\NN^r$-graded polynomial
$R$-algebra
$S$ and a compatible sequence of embeddings
$X_n\subset \Proj(S)$.  Choose an open neighborhood
$U\subset\Proj(S)$ so that $X_n\subset U$ are closed embeddings.
Let
$I_n\subset S$ be the  graded ideal of the closed embedding
$\overline{X}_n\subset \Proj(S)$, and set
$I=\bigcap_{n\geq 0} I_n$.
Then
$X=\Proj(S/I)\cap U$ is the desired algebraization of the formal scheme
$\foX$.
\qed

\begin{remark}
If we have $\shK=\shL^d$, then $P_0=X_0$, such that
the formal sheaf $\shL$ is ample on the closed fiber $X_0$.
In this case,
Grothendieck's Algebraization Theorem ensures
that $\foX$ is algebraizable.
\end{remark}

\begin{question}
The assumptions in Theorem \ref{algebraization} imply  that
the restriction of the family $\shL_1,\ldots,\shL_r$ to the closed fiber
is  ample. A natural question to ask:
Given
a proper formal scheme with
a family of invertible  sheaves whose restriction to
the closed fiber is ample -- is the formal scheme algebraizable?
\end{question}

%===========================================================

\end{document}